\documentclass[12pt,a4paper]{article}
\title{Convex Sets}

\usepackage[utf8]{inputenc}
\usepackage{amsmath}
\usepackage{amsfonts}
\usepackage{amssymb}
\usepackage{mathrsfs}
\usepackage{pbox}
\usepackage{accents}
\usepackage{url}
\usepackage{tabularx}
\usepackage[pdftex]{graphicx}
\usepackage[hang,flushmargin]{footmisc}
\usepackage{titlesec}
\usepackage[left=1.0in, right=1.0in, top=1.0in,bottom=1.0in]{geometry}
\usepackage[labelfont=bf,justification=centering,singlelinecheck=true,font=small]{caption}

\titleformat*{\subsection}{\normalsize\bfseries}
\titleformat{\section}{\large\bfseries}{\thesection}{1em}{}
\numberwithin{equation}{section}
\allowdisplaybreaks
\newtheorem{theorem}{Theorem}
\newtheorem{lemma}{Lemma}
\newtheorem{corollary}{Corollary}

\begin{document}

\begin{center}
\LARGE{Classification on convex sets in the presence of missing covariates}\\[14pt]
\Large
Levon Demirdjian\footnote{Corresponding author\\  Email: levondem@ucla.edu \\ Department of Statistics, University of California Los Angeles, USA 
\vspace{2pt}},   
Majid Mojirsheibani\footnote{Department of Mathematics, California State University Northridge, USA}
\normalsize
\end{center}

%~~~~~~~~~~~~~~~~~~~~~~~~~~~~~~~~~~~~~~~~~~~~~~~~~~~~%
%~~~~~~~~~~~~~~~~~~ ABSTRACT ~~~~~~~~~~~~~~~~~~~~~~~~%
%~~~~~~~~~~~~~~~~~~~~~~~~~~~~~~~~~~~~~~~~~~~~~~~~~~~~%

\hrulefill

\textbf{Abstract}\\
A number of results related to statistical classification on convex sets are presented. In particular, the
focus is on the case where some of the covariates in the data and observation being classified
can be missing. The form of the optimal classifier is derived when the class-conditional densities
are uniform over convex regions. In practice, the underlying convex sets are often unknown and
must be estimated with a set of data. In this case, the convex hull of a set of points is shown
to be a consistent estimator of the underlying convex set. The problem of estimation is further
complicated since the number of points in each convex hull is itself a random variable. The
corresponding plug-in version of the optimal classifier is derived and shown to be Bayes consistent.  \\
\


\hrulefill

%~~~~~~~~~~~~~~~~~~~~~~~~~~~~~~~~~~~~~~~~~~~~~~~~~~~~%
%~~~~~~~~~~~~~~~~~~ INTRODUCTION ~~~~~~~~~~~~~~~~~~~~%
%~~~~~~~~~~~~~~~~~~~~~~~~~~~~~~~~~~~~~~~~~~~~~~~~~~~~%

\section{Introduction}

Consider the following two-group classification problem. Let $(\textbf{X},Y)$ be a 
random pair with underlying distribution $F_{\textbf{X},Y}$, where the observed vector $\textbf{X} \in \mathbb{R}^{d}$ is used to predict the unknown class membership of $Y \in \{0,1 \}$. In classification, one seeks to find a function (classifier) $\phi :\mathbb{R}^{d} \to \{ 0,1\}$ whose probability of incorrect prediction, denoted
	\begin{equation*}
	L(\phi)=P\{\phi(\textbf{X}) \neq Y\},
	\end{equation*}
is as small as possible. The rule that minimizes this error is known as the 
\textit{Bayes classifier}, and is given by
	\begin{equation}
	\label{Bayes Rule}
	\phi_{B}(\textbf{x})=
	\begin{cases}
	1 \text{ if } P\{Y=1 \vert \textbf{X}=\textbf{x}\} > P\{Y=0 \vert \textbf{X}= 
	\textbf{x}\}\\
	0 \text{ otherwise, }
	\end{cases}
	\end{equation}
with corresponding probability of misclassification 
	\begin{equation*}
	L(\phi_{B})=P\{\phi_{B}(\textbf{X}) \neq Y\};
	\end{equation*}
that is, $L(\phi_{B}) \leq L(\phi)$ for all classifiers $\phi$. For more on the general problem of statistical classification, one can refer to Devroye et al. (1996) and their list of references.\\
\

There are several open questions related to classification when the class-conditional densities of $\textbf{X}$ given $Y$ are assumed to be known. For example, Devroye et al. (1996) examine classification when the class-conditional densities are uniform over the $d$-dimensional hyper-rectangles $\mathcal{R}_1$ and $\mathcal{R}_0$; that is, $\textbf{X}|Y=1 \sim$ Unif$(\mathcal{R}_{1})$ and $\textbf{ X}|Y=0 \sim$ Unif$(\mathcal{R}_{0}).$ If $\textbf{x}$ falls in $\mathcal{R}_{1} - \mathcal{R}_{0}$, for example, there is no probability of incorrect prediction; simply classify $\textbf{x}$ as belonging to class $1$. When $\textbf{x} \in \mathcal{R}_{1} \cap \mathcal{R}_{0}$, however, the problem becomes more challenging. In this setup, the authors show that the optimal rule in (\ref{Bayes Rule}) becomes
	
\begin{equation}
\label{Bayes Rule - Rectangles}
\phi_{B}(\textbf{x})=
\begin{cases}
1 \text{ if } \textbf{x} \in \mathcal{R}_1 - \mathcal{R}_0\\
0 \text{ if } \textbf{x} \in \mathcal{R}_0 - \mathcal{R}_1\\
1 \text{ if } \textbf{x} \in \mathcal{R}_1 \cap \mathcal{R}_0, \dfrac{p}{\mu(\mathcal{R}_1)} > \dfrac{1-p}{\mu(\mathcal{R}_0)}\\[10pt]
0 \text{ if } \textbf{x} \in \mathcal{R}_1 \cap \mathcal{R}_0, \dfrac{p}{\mu(\mathcal{R}_1)} \leq \dfrac{1-p}{\mu(\mathcal{R}_0)},\\[10pt]
\end{cases}
\end{equation}

where $p=P\{Y=1\}$ and $\mu(\cdot)$ denotes the Lebesgue measure. In practice, however, the assumption that the conditional densities are uniform over hyper-rectangles is far too restrictive and limits the use of the above classifier. Moreover, the regions $\mathcal{R}_1$ and $\mathcal{R}_0$ are almost always unknown and must therefore be estimated using a set of iid data. To complicate matters further, there is no mention of how to deal with the more realistic case where there may be missing covariates in both the data and new observation $\textbf{X}$ being classified. This paper aims to address all of these concerns.\\
\

Accurate estimation of the unknown regions $\mathcal{R}_1$ and $\mathcal{R}_0$ is far more than just a theoretical concern and has widespread applications in many diverse disciplines. In conservational biology, for example, Kerley et al. (2002) use location data collected from radio-collared tigers in parts of eastern Russia to estimate their home ranges. In a related work, Burgman and Fox (2003) propose using a generalization of the convex hull of observed locations to predict the shape of a species' true habitat. The cited works empirically verify that species' habitats almost always form complex geometrical regions. To conduct reliable inference, therefore, there is clearly a need to both generalize the underlying assumptions of the classifier in (\ref{Bayes Rule - Rectangles}) and find an optimal method of estimation.\\
\

In order to construct a sample-based version of the optimal rule in (\ref{Bayes Rule - Rectangles}) using estimates of the unknown regions $\mathcal{R}_1$ and $\mathcal{R}_0$, one typically uses an iid set of data  \linebreak $D_{n} = \{ (\textbf{X}_{1}, Y_{1}), ..., (\textbf{X}_{n}, Y_{n}) \}$, where $(\textbf{X}_{i}, Y_{i}) \sim F_{\textbf{X},Y}$ for $i=1,...,n$. Denoting the resulting data-based classifier as $\phi_{n}$, the hope is that the performance of $\phi_{n}$ is similar in some sense to that of the optimal classifier $\phi_{B}$. In general, the quality of a sample-based classifier is usually assessed by checking whether or not it is consistent; a data-based rule $\phi_n$ is said to be consistent (weakly) if  
	\begin{equation*}
	P\{\phi_n({ \bf X})\neq Y\,|\, D_n\} \overset{p} \longrightarrow P\{\phi_B({ 
	\bf X})\neq Y\}.
	\end{equation*}
The rule is said to be strongly consistent if the convergence holds with probability one (Devroye et al. (1996)).\\
\

In practice, the situation is more complicated when some of the covariates used to predict the class membership of $\textbf{X}$ are missing. Prior to 2007, most of the existing results in the literature on classification with missing covariates deal with the case where covariates can be missing in the data (i.e., in ${\bf X}_i,~i=1, ..., n$) but not in the new observation ${\bf X}$. See, for example, Chung and Han (2000) for the parametric case, and Pawlak (1993) for the nonparametric case. Here, however, we consider the more general case where there may be missing covariates in both the data and the new observation being classified. To see why missing values can cause difficulty, first note that if there are missing covariates in the data, the number of regions to be estimated (one corresponding to each class-conditional density) will typically increase (a more detailed explanation will be given in Section 2). Yet another source of difficulty is that in the case with missing covariates, the rule $\phi_{B}$ given in (\ref{Bayes Rule}) will not necessarily be the optimal classifier any more. In fact, Mojirsheibani and Montazeri (2007) derived the optimal classifier in the presence of missing covariates, which is in general different from (\ref{Bayes Rule}). For more recent results on classification with missing covariates, one may also refer to Mojirsheibani (2012) and Mojirsheibani and Chenouri (2011).\\

In this paper, we expand upon the above problem of classification when the class-conditional densities are known, and extend our results to the more realistic case where covariates can be missing in both the data and new observation being classified. In Section 2, we present the form of the classifier that is optimal not only in the restrictive case where $\mathcal{R}_{1}$ and $\mathcal{R}_{0}$ are hyper-rectangles, but in the more general situation where $\mathcal{R}_{1}$ and $\mathcal{R}_{0}$ can be any two convex sets. In Section 3, we propose a data-based version of the optimal rule which is constructed with sample-based approximations of the two convex regions. The proposed classifier is then shown to be strongly consistent. A formal justification for our method of estimating convex regions can be found in the appendix.\\
\

\section{Estimation}

Our discussion and results for classification on convex sets are based on the following setup. Let $\mathcal{C}\subset \mathbb{R}^d$ be convex and let $D_n=\{\textbf{X}_1,\textbf{X}_2,...,\textbf{X}_n\}$ denote an iid set of $d$-dimensional random vectors, where $\textbf{X}_i$ is uniformly distributed on $\mathcal{C}$ for $i=1,...,n$. Here and throughout the paper, we will be assuming that $\mathcal{C}$ is closed (this is not as restrictive as it seems; if $\mathcal{C}$ is not closed, one can instead consider the closure of $\mathcal{C}$, denoted $\bar{\mathcal{C}}$, since $P\{ \textbf{X} \in \mathcal{C} \} = P\{ \textbf{X} \in \bar{\mathcal{C}} \}$). Given $D_n$, how can we estimate $\mathcal{C}$ in a consistent way? To motivate our choice of using the convex hull of $D_{n}$ as an estimator of $\mathcal{C}$, we consider the method of maximum likelihood. Let $\theta$ be the (unknown) volume of $\mathcal{C} $ and note that the likelihood function can be written as

	\begin{equation*}
	L(\theta \vert D_n)=\prod_{i=1}^{n}\dfrac{I\{\textbf{x}_i\in 					\mathcal{C} \}}{\theta} =\dfrac{I\{ \textbf{x}_1,\textbf{x}_2,...,				\textbf{x}_n \in \mathcal{C}\}}{{\theta}^n},
	\end{equation*}

where $I\{ \cdot \}$ denotes the indicator function ($I\{ A \}=1$ if $\textbf{x} \in A$ and $I\{ A\} = 0$ if $\textbf{x} \not\in A$ for any set $A \subset \mathbb{R}^{d}$). Clearly, this expression is maximized when our estimate of $\theta$ is minimized. With this in mind we take the convex hull of $D_n$, which we shall denote by $\widehat{\mathcal{C}}_n$, to be our estimate of $\mathcal{C}$ (the convex hull of $D_n$ is, by definition, the smallest convex set containing the points in $D_n$).  How good is $\widehat{\mathcal{C}}_n$ at estimating $\mathcal{C}$? First, let us show that $\widehat{\mathcal{C}}_n \overset{a.s.} \longrightarrow \mathcal{C}$ in some sense.\\
	
For a fixed point $\textbf{t}\in \mathcal{C}$, define the ``distance" between $\textbf{t}$ and $\widehat{\mathcal{C}}_n$ to be $ \inf_{y\in \widehat{\mathcal{C}}_n}\|\textbf{t}-\textbf{y}\|,$ where $\| \cdot \|$ denotes the Euclidean norm in $\mathbb{R}^d$. To establish the convergence of $\widehat{\mathcal{C}}_n$ to $\mathcal{C}$, we will consider the Hausdorff distance 

\begin{equation*}
d(\widehat{\mathcal{C}}_n, \mathcal{C}) = \max\left\lbrace                    \sup\limits_{\textbf{t} \in \mathcal{C}} \inf\limits_{y\in \widehat{\mathcal{C}}_n}\|\textbf{t}-\textbf{y}\|, \sup\limits_{y\in \widehat{\mathcal{C}}_n} \inf\limits_{\textbf{t} \in \mathcal{C}}\|\textbf{t}-\textbf{y}\|\right\rbrace.
\end{equation*}

Since $\widehat{\mathcal{C}}_n \subset \mathcal{C}$, however, $d(\widehat{\mathcal{C}}_n, \mathcal{C})$ can be written equivalently as $\sup_{\textbf{t} \in \mathcal{C}} \inf_{y\in \widehat{\mathcal{C}}_n}\|\textbf{t}-\textbf{y}\|$. Therefore, to show that  $\widehat{\mathcal{C}}_n \overset{a.s.} \longrightarrow \mathcal{C}$, we will show that $\sup\limits_{\textbf{t} \in \mathcal{C}} \inf\limits_{y\in \widehat{\mathcal{C}}_n}\|\textbf{t}-\textbf{y}\| \overset{a.s.} \longrightarrow 0$ as $n \to \infty$.

	\begin{lemma}
	\label{Lemma - d(datapoint, t)}
	Let $\mathcal{C}\subset \mathbb{R}^d$ be convex and let 						$D_n=\{\textbf{X}_1,\textbf{X}_2,...,\textbf{X}_n\}$ denote a set of 			iid random vectors, uniformly distributed on $\mathcal{C}$. Also, denote 		the convex hull of $D_n$ by $\widehat{\mathcal{C}}_n$. Then for all 			$\textbf{t} \in \mathcal{C}$ and $\epsilon>0$, 
		\begin{equation*}
		P\left\lbrace \min_{i=1}^{n} \| \textbf{t}-\textbf{X}_i\|						>\epsilon\right\rbrace = r^{n}
		\end{equation*}
 	for some $r:= r(\epsilon) \in [0,1)$, where $r$ does not depend on $n$.
 	\end{lemma}
 
\textit{Proof of Lemma \ref{Lemma - d(datapoint, t)}}\\

Let $\textbf{t} \in \mathcal{C}$ and $\epsilon>0$ be given. Denote the $d-$dimensional ball with center $\textbf{t}$ and radius $\epsilon$ by $B(\textbf{t},\epsilon)$. If $\mu(A)$ represents the Lebesgue measure of a set $A$, we get

\begin{align*}
P\left\lbrace \min_{i=1}^{n} \| \textbf{t}-\textbf{X}_i\|>\epsilon\right\rbrace &= P\{\textbf{X}_1,...,\textbf{X}_n\not \in B(\textbf{t},\epsilon)\} = \prod_{i=1}^{n}P\{ \textbf{X}_i \not \in B(\textbf{t},\epsilon)\}\\
&=\left[ P\{ \textbf{X}_1 \not \in B(\textbf{t},\epsilon)\} \right]^{n} =\left[ 1-\dfrac{\mu[B(\textbf{t},\epsilon)\cap \mathcal{C} ]}{\mu(\mathcal{C})} \right]^{n}\\
&:=r^{n}.
\end{align*}

Since $\textbf{t} \in \mathcal{C}$, $\mu[B(\textbf{t},\epsilon)\cap \mathcal{C} ] > 0$, and so $0 \leq r < 1. \hspace{15pt} \Box$\\
\

\begin{lemma}
\label{Lemma - d(datapoint, t) RANDOM}
Let $\mathcal{C}\subset \mathbb{R}^d$ be convex and let $N \sim \text{binomial}(n,p)$. For all $k \in  \{1, 2,... , n\}$, let $D_k=\{\textbf{X}_1,\textbf{X}_2,...,\textbf{X}_k\}$ be a set of iid random vectors uniformly distributed on $\mathcal{C}$ and denote the convex hull of $D_k$ by $\widehat{\mathcal{C}}_k$. Then for all $\textbf{t} \in \mathcal{C}$ and $\epsilon > 0$, 

\begin{equation*}
\sum_{k=1}^{n}P\left\lbrace \inf_{ \textbf{y} \in \widehat{\mathcal{C}}_ k} \| \textbf{t}-\textbf{y}\| > \epsilon \right\rbrace P \{ N = k \} \leq c_{1}^{n} + c_{2}^{n},
\end{equation*}
 
where $c_{1}: = c_{1}(\epsilon) \in (0,1)$, $c_{2}: = c_{2}(\epsilon) \in (0,1)$, and neither $c_{1}$ nor $c_{2}$ depends on $n$. 
 
\end{lemma}

\textit{Proof of Lemma \ref{Lemma - d(datapoint, t) RANDOM}}\\
 
Let $\textbf{t} \in \mathcal{C}$ and $\epsilon>0$ be given. Then

\begin{align*}
\sum_{k=1}^{n} P\left\lbrace \inf_{ \textbf{y} \in \widehat{\mathcal{C}}_k} \| 	\textbf{t}-\textbf{y}\| > \epsilon \right\rbrace P\{ N=k \}
&\leq \sum_{k=1}^{n} P\left\lbrace \min_{i=1}^{k}\| \textbf{t}-\textbf{X}_i\| > \epsilon  \right\rbrace P\{ N=k \}\\
&= \sum_{k=1}^{n} r^k {n\choose k} p^k (1-p)^{n-k} \hspace{10pt} \left (\text{$0 \leq r < 1$; by Lemma \ref{Lemma - d(datapoint, t)}} \right)\\
&=(1-p)^n \sum_{k=1}^{n} {n\choose k} \left( \frac{rp}{1-p} \right)^k 1^{n-k}\\
&=(1-p)^n \left[ \left( 1+ \frac{rp}{1-p}\right)^n - 1\right]\\
&=[1-p+rp]^n - (1-p)^n\\
&:=c_{1}^{n} + c_{2}^{n}, 
\end{align*}
 
where $0 < c_{1} < 1$, $0 < c_{2} < 1$, and neither $c_{1}$ nor $c_{2}$ depends on $n$.  \hspace{15pt} $\Box$

\begin{lemma}
\label{Lemma - compact}
If $\mathcal{C}\subset \mathbb{R}^d$ is convex and closed with finite Lebesgue measure, then $\mathcal{C}$ is compact.
\end{lemma}

\textit{Proof of Lemma \ref{Lemma - compact}}\\

Since $\mathcal{C}$ is convex and has probability measure $1$, it is bounded.
By assumption, $\mathcal{C}$ is also closed, and therefore compact.\hspace{15pt} $\Box$\\

Up to this point, we have been assuming that the size of our sample $D_{n}$ was fixed. When we discuss classification in the next section, however, we will need to consider the more general case where the sample size is itself a random variable. To see why, note that we can always partition our data into two sets - those $\textbf{X}_{i}'$s for which $Y_{i}=1$ and those $\textbf{X}_{i}'$s for which $Y_{i}=0$. Since there is no way of knowing how many points will lie in each set, the sample sizes of both sets will be random.\\
\

Before stating our main result, we must first define the convex hull of a set of points containing a random number of elements. Let $N$ be a non-negative integer valued random variable and let $D_{N} = \{ \textbf{X}_{1}, \textbf{X}_{2}, ... , \textbf{X}_{N}\}.$ If $N>0$, define $\widehat{C}_{N}$ to be the convex hull of $D_{N}$; if $N=0$, define $\widehat{C}_{N}$ to be a ball of arbitrary radius $\delta > 0$ about the origin.\\
\

\begin{theorem}
\label{Theorem - Hull to True convergence}
Let $\mathcal{C}\subset \mathbb{R}^d$ be convex and closed and let $N \sim \text{binomial}(n,p)$. Define $D_N=\{\textbf{X}_1,\textbf{X}_2,...,\textbf{X}_N\}$ to be a set of iid random vectors, uniformly distributed on $\mathcal{C}$. Also, denote the convex hull of $D_N$ by $\widehat{\mathcal{C} }_N$, as defined above. Then
\begin{equation*}
\widehat{\mathcal{C}}_N \overset{a.s.} \longrightarrow \mathcal{C}
\end{equation*}
in the Hausdorff metric.
\end{theorem}
\textit{Proof of Theorem \ref{Theorem - Hull to True convergence}}\\

Let $\epsilon > 0$ and observe that

	\begin{align}
	\label{Theorem 1 - Piece 0}
	P \left\lbrace \sup_{\textbf{t} \in \mathcal{C}} 						\inf\limits_{\textbf{y} 		\in \widehat{\mathcal{C}}_N}\|			\textbf{t}-\textbf{y}\| > \epsilon 				\right\rbrace &= P 		\left\lbrace \sup_{\textbf{t} \in \mathcal{C}} 							\inf\limits_{\textbf{y} \in \widehat{\mathcal{C}}_N}\|\textbf{t} - 		\textbf{y} \| > \epsilon, N>0 \right\rbrace \notag \\
	&+ P \left\lbrace \sup_{\textbf{t} \in \mathcal{C}} \inf\limits_{ 		\textbf{y} \in \widehat{\mathcal{C}}_N}\|\textbf{t} - \textbf{y} \| 	> \epsilon, N=0 \right\rbrace \notag \\
	& \leq P \left\lbrace \sup_{\textbf{t} \in \mathcal{C}} 				\inf\limits_{\textbf{y} \in \widehat{\mathcal{C}}_N}\|\textbf{t} - 		\textbf{y} \| > \epsilon, N>0 \right\rbrace  + P\{ N=0 \} \notag \\
	& = \sum_{\ell=1}^{n} P \left\lbrace \sup_{\textbf{t} \in 				\mathcal{C}} \inf\limits_{\textbf{y} \in \widehat{\mathcal{C} 			}_{N}}\|\textbf{t} - \textbf{y}\| > \epsilon \Big| N= \ell 				\right\rbrace P \{ N = \ell\} + P\{ N=0\} \notag \\
	& = \sum_{\ell=1}^{n} P_{\ell, \mathcal{C}} \left\lbrace 				\sup_{\textbf{t} \in \mathcal{C}} \inf\limits_{\textbf{y} \in 			\widehat{\mathcal{C} }_{\ell}}\|\textbf{t} - \textbf{y}\| > 			\epsilon  \right\rbrace P \{ N = \ell\} + P\{ N=0\} \notag \\
	& := \text{I + II},
	\end{align}

where $P_{\ell, \mathcal{C}}\{ A \}$ denotes the probability of the event $A$ under the joint uniform distribution of $(\textbf{X}_{1}, \textbf{X}_{2}, ... ,\textbf{X}_{\ell})$ on $\mathcal{C}$. To deal with term I, consider the collection $\{B(\textbf{t},\frac{\epsilon}{2}) \vert \textbf{t} \in \mathcal{C}\}$, i.e. the set of all open balls of radius $\frac{\epsilon}{2}$ with centers in $\mathcal{C}$. Clearly, $\mathcal{C} \subset \cup_{\textbf{t} \in \mathcal{C}}B(\textbf{t}, \frac{\epsilon}{2} )$. Furthermore, since $\mathcal{C}$ is compact (see Lemma \ref{Lemma - compact}), there exists a finite subcover of $\mathcal{C}$, say $\cup_{i=1}^{k}B(\textbf{t}_{i},\frac{\epsilon}{2})$. For any $\ell \in \{ 1,2,...,n\}$,
	
	\begin{equation}
	\label{Theorem 1 - Piece 1}
	\sup\limits_{\textbf{t} \in \mathcal{C}} \inf\limits_{\textbf{y}\in 	\widehat{\mathcal{C}}_\ell}\|\textbf{t}-\textbf{y}\|\leq 				\max\limits_{i=1 }^{k} \sup\limits_{\textbf{t} \in B(\textbf{t}_i,		\frac{\epsilon}{2})} \inf\limits_{\textbf{y}\in \widehat{ 				\mathcal{C} }_\ell}\|\textbf{t} -\textbf{y}\|.
	\end{equation} 

By the triangle inequality, 

\begin{equation*}
\sup\limits_{\textbf{t} \in B(\textbf{t}_i,\frac{\epsilon}{2})} \inf\limits_{\textbf{y}\in \widehat{\mathcal{C}}_\ell}\|\textbf{t}-\textbf{y}\|\leq \frac{\epsilon}{2} + \inf\limits_{\textbf{y}\in \widehat{\mathcal{C}}_\ell}\|\textbf{t}_i-\textbf{y}\|
\end{equation*}

for all $i\in\{ 1,2,...,k\}$. Therefore,
\begin{align}
\label{Theorem 1 - Piece 2}
\max_{i=1}^{k}\sup\limits_{\textbf{t} \in B(\textbf{t}_i,\frac{\epsilon}{2})} \inf\limits_{\textbf{y}\in \widehat{\mathcal{C}}_\ell}\|\textbf{t}-\textbf{y}\| &\leq \frac{\epsilon}{2} + \max_{i=1}^{k} \inf\limits_{\textbf{y}\in \widehat{\mathcal{C}}_\ell}\|\textbf{t}_i-\textbf{y}\| \leq \frac{\epsilon}{2} + \sum_{i=1}^{k} \inf\limits_{\textbf{y}\in \widehat{\mathcal{C}}_\ell}\|\textbf{t}_i-\textbf{y}\|.
\end{align}

Putting expressions (\ref{Theorem 1 - Piece 1}) and (\ref{Theorem 1 - Piece 2}) together yields  
	\begin{align*}
	P_{\ell, \mathcal{C}} \left\lbrace \sup_{\textbf{t} \in \mathcal{C}} \inf\limits_{ \textbf{y} 	\in \widehat{\mathcal{C}}_\ell}\|\textbf{t} - \textbf{y}\| > \epsilon 			\right\rbrace &\leq P_{\ell, \mathcal{C}} \left\lbrace \frac{\epsilon}{2} + \sum_{i=1}^{k} 			\inf\limits_{\textbf{y}\in \widehat{\mathcal{C}}_\ell}\|\textbf{t}_i - 			\textbf{y}\| > \epsilon \right\rbrace\\
	& \leq \sum_{i=1}^{k} P_{\ell, \mathcal{C}} \left\lbrace \inf\limits_{ \textbf{y}\in 				\widehat{\mathcal{C}}_\ell}\|\textbf{t}_i - \textbf{y}\| > 						\frac{\epsilon}{2k} \right\rbrace
	\end{align*}

for all $\ell \in\{ 1,2,...,n\}$. We can now bound term I in expression (\ref{Theorem 1 - Piece 0}) via
	\begin{align*}
	\sum_{\ell=1}^{n} P_{\ell, \mathcal{C}} \left\lbrace \sup_{\textbf{t} \in \mathcal{C}} 				\inf\limits_{\textbf{y} \in \widehat{\mathcal{C}}_\ell}\|\textbf{t} - 			\textbf{y}\| > \epsilon \right\rbrace P \{ N = \ell\} & \leq \sum_{i = 	    	1}^{k} \sum_{\ell=1}^{n} P_{\ell, \mathcal{C}} \left\lbrace \inf\limits_{ \textbf{y}\in 			\widehat{\mathcal{C}}_\ell}\|\textbf{t}_i - \textbf{y}\| > 						\frac{\epsilon}{2k} \right\rbrace P \{ N = \ell\}\\
	& \leq \sum_{i=1}^{k} (c_{1}^{n} + c_{2}^{n}) \hspace{15pt} \text{(by 			Lemma \ref{Lemma - d(datapoint, t) RANDOM})}
	\end{align*}

where $0 < c_{1} < 1$, $0 < c_{2} < 1$, and neither $c_{1}$ nor $c_{2}$ depends on $n$. To deal with term II in expression (\ref{Theorem 1 - Piece 0}), simply note that 

\begin{equation*}
P \{ N=0 \} = (1-p)^{n} = e^{n\log(1-p)} = e^{-c_{3}n},
\end{equation*}

where $c_{3} > 0$. Putting all of the above together, we have

\begin{equation*}
P \left\lbrace \sup_{\textbf{t} \in \mathcal{C}} \inf\limits_{\textbf{y} 		\in \widehat{\mathcal{C}}_N}\|\textbf{t}-\textbf{y}\| > \epsilon 				\right\rbrace \leq kc_{1}^{n} + kc_{2}^{n} + e^{-c_{3}n}.
\end{equation*}

Therefore

\begin{equation*}
\sum_{n=1}^{\infty} P \left\lbrace \sup_{\textbf{t} \in \mathcal{C}} \inf\limits_{\textbf{y} \in \widehat{\mathcal{C}}_N}\|\textbf{t}-\textbf{y}\| > \epsilon 	\right\rbrace \leq k \sum_{n=1}^{\infty} c_{1}^{n} + k \sum_{n=1}^{\infty} c_{2}^{n} + \sum_{n=1}^{\infty} e^{-c_{3}n} < \infty,
\end{equation*}

and by the Borel-Cantelli lemma, $\sup\limits_{\textbf{t} \in \mathcal{C}} \inf\limits_{\textbf{y} \in \widehat{\mathcal{C}}_N}\|\textbf{t}-\textbf{y}\| \overset{a.s.} \longrightarrow 0$ as $n \to \infty$. $\hspace{15pt} \Box$

%%%%%%%%%%%%%%%%%%%%%%%%%%%%%%%%%%%%%%%%%%%%%%%%%%%%

\section{Classification on Convex Sets}

%%%%%%%%%%%%%%%%%%%%%%%%%%%%%%%%%%%%%%%%%%%%%%%%%%%%

Now that we have shown that $\widehat{\mathcal{C}}_{n}$ is a consistent estimator of $\mathcal{C}$, we can carefully examine classification on convex sets.  Before considering our main results, i.e. those involving classification with missing covariates, we deal with the simpler case where every observation is available. The ideas in this section will make our main results in the following section more presentable.\\
\

Let $\textbf{X} \in \mathbb{R}^{d}$ be a vector of covariates to be used to predict the class membership $Y\in \{0,1\}$. Here, $\textbf{X} \vert Y=1 \sim $ Unif $(\mathcal{C}_1)$ and $\textbf{X} \vert Y=0 \sim $ Unif $(\mathcal{C}_0)$, where $\mathcal{C}_1$ and $\mathcal{C}_0 $ are taken to be convex sets in $\mathbb{R}^d$. In this setup, the optimal rule $\phi_{B}$ is given by (\ref{Bayes Rule - Rectangles}), with $\mathcal{R}_1$ and $\mathcal{R}_0$ replaced by $\mathcal{C}_1$ and $\mathcal{C}_0$ respectively. To see this, rewrite expression (\ref{Bayes Rule}) as 
\begin{align*}
	\phi_{B}(\textbf{x})&= I \left\lbrace P\{Y=1\vert \textbf{X} = 					\textbf{x}\} >  P\{Y=0 \vert \textbf{X} = \textbf{x}\} \right\rbrace\\
	&=I \left\lbrace \textbf{x} \in \mathcal{C}_{1} \cap \mathcal{C}_{0}, 			P\{Y=1\vert \textbf{X} = \textbf{x}\} > P\{Y=0 \vert \textbf{X} = 				\textbf{x}\} \right\rbrace\\
	&+I \left\lbrace \textbf{x} \not\in \mathcal{C}_{1} \cap \mathcal{C}_{0}, 	P\{Y=1\vert \textbf{X} = \textbf{x}\} > P\{Y=0 \vert \textbf{X} = 				\textbf{x}\} \right\rbrace\\
	&:= \text{I + II}.
	\end{align*}

Expand the second term to get 
\begin{align*}
	\text{II}&= I \left\lbrace \textbf{x} \in \mathcal{C}_{1} - 					\mathcal{C}_{0}, P\{Y=1\vert \textbf{X} = \textbf{x}\} > P\{Y=0 \vert 			\textbf{X} = \textbf{x}\} \right\rbrace\\
	&+  I \left\lbrace \textbf{x} \in \mathcal{C}_{0} - \mathcal{C}_{1}, 			P\{Y=1\vert \textbf{X} = \textbf{x}\} > P\{Y=0 \vert \textbf{X} = 				\textbf{x}\} \right\rbrace\\
	&=I \left\lbrace \textbf{x} \in \mathcal{C}_{1} - 								\mathcal{C}_{0}\right\rbrace + 0\\
	&=I \left\lbrace \textbf{x} \in \mathcal{C}_{1} - 								\mathcal{C}_{0}\right\rbrace.
	\end{align*}

Therefore we see that
\begin{align}
\label{Bayes Rule Indicator - Convex (no missing)}
\phi_{B}(\textbf{x})&= I \left\lbrace \textbf{x} \in \mathcal{C}_{1} \cap 	\mathcal{C}_{0}, P\{Y=1\vert \textbf{X} = \textbf{x}\} > P\{Y=0 \vert 			\textbf{X} = \textbf{x}\} \right\rbrace \nonumber \\ 
&+ I \left\lbrace \textbf{x} \in \mathcal{C}_{1} - 								\mathcal{C}_{0}\right\rbrace.
\end{align}

Denote the probability density function of $\textbf{X}$ conditioned on $Y=y$ by $f_{\textbf{X}\vert Y=y}$, and let $f_\textbf{X}$ be the marginal density function of $\textbf{X}$. Next, let $p=P\{Y=1\}$ and let $\mu(\cdot)$ denote the Lebesgue measure of a set and use the fact that whenever $\textbf{x} \in \mathcal{C}_{1} \cap \mathcal{C}_{0}$,

\begin{equation}
	\label{Conditional 1}
	P\{Y=1 \vert \textbf{X}=\textbf{x}\}=\dfrac{f_{\textbf{X}\vert Y=1}				(\textbf{x} \vert 1) P\{Y=1\}}{f_\textbf{X}(\textbf{x})} = 						\dfrac{\dfrac{p} {\mu{(\mathcal{C}_{1})}}}{\dfrac{p}							{\mu{(\mathcal{C}_{1})}} + \dfrac{1-p}	{\mu{(\mathcal{C}_{0})}}}
	\end{equation}

and

	\begin{equation}
	\label{Conditional 2}
	P\{Y=0 \vert \textbf{X}=\textbf{x}\}=\dfrac{f_{\textbf{X}\vert Y=0}				(\textbf{x} \vert 0) P\{Y=0\}}{f_\textbf{X}(\textbf{x})} = 						\dfrac{\dfrac{(1-p)}{\mu{(\mathcal{C}_{0})}}}{\dfrac{p}							{\mu{(\mathcal{C}_{1})}} + \dfrac{1-p}{\mu{(\mathcal{C}_{0})}}}
	\end{equation}

to rewrite (\ref{Bayes Rule Indicator - Convex (no missing)}) as

\begin{equation}
\label{Bayes Rule - Convex}
\phi_{B}(\textbf{x})=
\begin{cases}
1 \text{ if } \textbf{x} \in \mathcal{C}_1 - \mathcal{C}_0\\
0 \text{ if } \textbf{x} \in \mathcal{C}_0 - \mathcal{C}_1\\
1 \text{ if } \textbf{x} \in \mathcal{C}_1 \cap \mathcal{C}_0, \dfrac{p}{\mu(\mathcal{C}_1)} > \dfrac{1-p}{\mu(\mathcal{C}_0)}\\[10pt]
0 \text{ if } \textbf{x} \in \mathcal{C}_1 \cap \mathcal{C}_0, \dfrac{p}{\mu(\mathcal{C}_1)} \leq \dfrac{1-p}{\mu(\mathcal{C}_0)}.\\[10pt]
\end{cases}
\end{equation}

Thus, we have proved the following result.

\begin{theorem}
\label{Theorem - Optimal Classifier}
Let $\textbf{X} \in \mathbb{R}^{d}$ be a vector of covariates to be used to predict the class membership $Y\in \{0,1\}$, where $\textbf{X} \vert Y=1 \sim $ Unif $(\mathcal{C}_1)$ and $\textbf{X} \vert Y=0 \sim $ Unif $(\mathcal{C}_0)$ for convex sets $\mathcal{C}_1$ and $\mathcal{C}_0 $ in $\mathbb{R}^d$. Then the classifier which minimizes the probability of incorrect prediction is given by $\phi_{B}$ in (\ref{Bayes Rule - Convex}).
\end{theorem}

Since $\mathcal{C}_1$ and $\mathcal{C}_0$ are generally unknown, we may consider a data-based version of the classifier above. Let $D_n=\{(\textbf{X}_1,Y_1),...,(\textbf{X}_n,Y_n)\}$ be a set of iid data, distributed as $(\textbf{X},Y)$, where $\textbf{X}_i \in \mathbb{R}^d$ and $Y_i \in \{0,1\}$ for $i=1,...,n$. Here, $\textbf{X}_i \vert Y_i=1 \sim$Unif$(\mathcal{C}_1)$ and $\textbf{X}_i \vert Y_i = 0\sim$Unif$(\mathcal{C}_0)$ for $i=1,...,n$. Next, let $N_{1} \sim \text{binomial}(n,p)$ and $N_{0} = (n - N_{1}) \sim \text{binomial}(n,1-p)$ be the number of class $1$ and class $0$ points respectively (i.e. $N_{k}$ is the number of $(\textbf{X}_{i}, Y_{i})$ pairs for which $Y_{i}=k$, $k=0,1$). Also, let $\widehat{\mathcal{C}}_{1,N_{1}}:= \widehat{\mathcal{C}}_1 $ denote the convex hull of the class $1$ points, i.e. $\widehat{\mathcal{C}}_1$ is the convex hull of those $\textbf{X}_{i}'$s for which $Y_{i} = 1$. Similarly, let $\widehat{\mathcal{C}}_{0,N_{0}}:= \widehat{\mathcal{C}}_0$ denote the convex hull of those $\textbf{X}_{i}'$s for which $Y_{i} = 0$. Consider the plug in version of the rule in (\ref{Bayes Rule - Convex}):

\begin{equation}
	\label{Bayes Convex Sets Plug In}
	\phi_{n}(\textbf{x})=
		\begin{cases}
		1 \text{ if } \textbf{x} \in \widehat{\mathcal{C}}_1 - 							\widehat{\mathcal{C}}_0\\
		0 \text{ if } \textbf{x} \in \widehat{\mathcal{C}}_0 - 							\widehat{\mathcal{C}}_1\\
		1 \text{ if } \textbf{x} \in \widehat{\mathcal{C}}_1 \cap 						\widehat{\mathcal{C}}_0, \dfrac{\widehat{p}}									{\mu(\widehat{\mathcal{C}}_1)} > \dfrac{1-\widehat{p}}							{\mu(\widehat{\mathcal{C}}_0)}\\[15pt]
		0 \text{ if } \textbf{x} \in \widehat{\mathcal{C}}_1 \cap 						\widehat{\mathcal{C}}_0, \dfrac{\widehat{p}}									{\mu(\widehat{\mathcal{C}}_1)} \leq \dfrac{1-\widehat{p}}						{\mu(\widehat{\mathcal{C}}_0)}\\[15pt]
		I\left\lbrace \widehat{p} \cdot \inf\limits_{\textbf{y}\in 						\widehat{\mathcal{C}}_1}\| \textbf{x}-\textbf{y}\| < 							(1-\widehat{p}) \cdot \inf\limits_{\textbf{y}\in 								\widehat{\mathcal{C}}_0} \|	\textbf{x}-\textbf{y}\|								\right\rbrace\hspace{15pt}\text{ if } \hspace{5pt} \textbf{x} \not\in 		\widehat{\mathcal{C}}_1 \cup \widehat{\mathcal{C} }_0,
		\end{cases}
	\end{equation}

where $\widehat{p}=\frac{1}{n}\sum_{i=1}^{n}Y_{i}$.\\
\

How good is the classifier $\phi_{n}$ for predicting $Y$? To answer this question, let 

	\begin{equation*}
	L_{n}(\phi_{n}) = P \{ \phi_{n}(\textbf{X}) \neq Y \vert D_{n} \} 		\hspace{15pt} \text{ and } \hspace{15pt} L(\phi_{B}) = P \{ 			\phi_{B} (\textbf{X}) \neq Y\}
	\end{equation*}

be the misclassification probabilities of $\phi_{n}$ and $\phi_{B}$ respectively. The following theorem establishes the strong consistency of the plug in rule $\phi_{n}$.

	\begin{theorem}
	\label{Consistency}
	Let $\phi_{n}$ be the classifier in equation (\ref{Bayes Convex 		Sets Plug 	In}). Then 	$L_{n}(\phi_{n}) \overset{a.s.} 				\longrightarrow L(\phi_{B})$ as $n \to 	\infty$, where $\phi_{B}$ 		is the optimal rule in (\ref{Bayes Rule - Convex}).
	\end{theorem}

A few lemmas will provide the tools that are necessary to prove Theorem \ref{Consistency}.

\begin{lemma} 
\label{Beer's lemma}
Let $\{ \widehat{\mathcal{C}}_{N} \}, \mathcal{C}$ be as in Theorem \ref{Theorem - Hull to True convergence}. Then
	\begin{equation*}
	\mu(\widehat{\mathcal{C}}_{N} \Delta \mathcal{C}) \overset{a.s.} 		\longrightarrow 0
	\end{equation*}
as $n \to \infty$, where $\widehat{\mathcal{C}}_{N} \Delta \mathcal{C}$ denotes the symmetric difference of $\widehat{\mathcal{C}}_{N}$ and $\mathcal{C}$ and where $\mu(\cdot)$ denotes the Lebesgue measure.
\end{lemma}

\textit{Proof of Lemma \ref{Beer's lemma}}\\
\

By Theorem \ref{Theorem - Hull to True convergence}, $\widehat{ \mathcal{C}}_{N} \overset{a.s.} \longrightarrow \mathcal{C}$ in the Hausdorff metric. Furthermore, a result of Beer (1974) implies that for each $\omega$ in the underlying sample space $\Omega$, if $\widehat{\mathcal{C}}_{N}(\omega) \to \mathcal{C}$ (in the Hausdorff metric), then $\mu(\widehat{\mathcal{C}}_{N}(\omega) \Delta \mathcal{C}) \to 0$. Let $\omega \in \left\lbrace \omega \in \Omega \Big| \widehat{\mathcal{ C}}_{N}(\omega) \to \mathcal{C} \right\rbrace$. Then  $$\omega \in \left\lbrace \omega \in \Omega \Big| \mu(\widehat{\mathcal{C}}_{N}(\omega) \Delta \mathcal{C}) \to 0  \right\rbrace. $$ The lemma now follows since $P\left\lbrace \omega \in \Omega \Big| \widehat{\mathcal{ C}}_{N}(\omega) \to \mathcal{C} \right\rbrace = 1. \hspace{15pt} \Box$

\begin{corollary}	
\label{Corollary - Symmetric Difference}
Let $\{ \widehat{\mathcal{C}}_{N} \}, \mathcal{C}$ be as in Theorem \ref{Theorem - Hull to True convergence}. Then	
	\begin{equation*}
	P \left\lbrace \textbf{X} \in \widehat{\mathcal{C}}_{N} \Delta 			\mathcal{C} \Big| D_{n} \right\rbrace \overset{a.s.} 					\longrightarrow 0
	\end{equation*}
as $n \to \infty$.	
\end{corollary}	
	
\textit{Proof of Corollary \ref{Corollary - Symmetric Difference}}

	\begin{align}
	\label{Corollary 1}
	P \left\lbrace \textbf{X} \in \widehat{\mathcal{C}}_{N} \Delta 			\mathcal{C} \Big| D_{n} \right\rbrace &= P \left\lbrace \textbf{X} 		\in \widehat{\mathcal{C}}_{N} \Delta \mathcal{C}, N > 0 \Big| D_{n} 	\right\rbrace + P \left\lbrace \textbf{X} \in \widehat{\mathcal{C 		}}_{N} \Delta \mathcal{C}, N=0 	\Big| D_{n} \right\rbrace \notag \\
	& \leq I \{ N>0 \} P \left\lbrace \textbf{X} \in \widehat{ 				\mathcal{C }}_{N} \Delta \mathcal{C} \Big| D_{n} \right\rbrace + I 		\{ N = 0 \}.  
	\end{align}

Since $\widehat{\mathcal{C}}_{N} \subset \mathcal{C}$ for $N>0$ and since $\textbf{X}$ is uniformly distributed on $\mathcal{C}$, the inequality in (\ref{Corollary 1}) becomes
	\begin{align*}
	P \left\lbrace \textbf{X} \in \widehat{\mathcal{C}}_{N} \Delta 			\mathcal{C} \Big| D_{n} \right\rbrace \leq I \{ N>0 \} 					\frac{\mu(\widehat{\mathcal{C}}_{N} \Delta \mathcal{C})}				{\mu(\mathcal{C})} + I \{ N = 0 \}.  
	\end{align*}

By Lemma \ref{Beer's lemma}, $\frac{\mu(\widehat{\mathcal{C}}_{N} \Delta \mathcal{C})}{\mu(\mathcal{C})} \overset{a.s.} \longrightarrow 0$ as $n \to \infty$. Finally, $I\{N=0\} \overset{a.s.} \longrightarrow 0$ as $n \to \infty$ by Markov's inequality in conjunction with the Borel-Cantelli lemma. $\hspace{15pt} \Box$

\begin{lemma} 
\label{Lemma - Symmetric Difference}
Let $\{ \widehat{\mathcal{C}}_{N} \}, \mathcal{C}$ be as in Theorem \ref{Theorem - Hull to True convergence}. Then
	\begin{equation*}
	\mu(\widehat{\mathcal{C}}_{N}) \overset{a.s.} \longrightarrow 			\mu(\mathcal{C})
	\end{equation*}
as $n \to \infty$.
\end{lemma}

\textit{Proof of Lemma \ref{Lemma - Symmetric Difference}}\\
	\begin{align*}
	\Big| \mu(\mathcal{C}) - \mu(\widehat{\mathcal{C}}_{N}) \Big|  &= 
	\Big| \mu(\mathcal{C}) - \mu(\widehat{\mathcal{C}}_{N}) \Big| I \{ 		N > 0 \} + \Big| \mu(\mathcal{C}) - \kappa_{0} \Big| I \{ N=0 \}, 
	\end{align*}
	
where $\kappa_{0} = \mu(B(\textbf{0}, \delta))$ for some $ 0 < \delta < \infty$ (recall how we defined $\widehat{\mathcal{C}}_{N}$ for $N=0$; see the definition before the statement of Theorem \ref{Theorem - Hull to True convergence}). Therefore,

	\begin{equation*}
	\Big| \mu(\mathcal{C}) - \mu(\widehat{\mathcal{C}}_{N}) \Big| \leq 		\left[ \mu(\widehat{\mathcal{C}}_{N} \Delta \mathcal{C}) + \Big| \mu(\mathcal{C}) - \kappa_{0} \Big| I \{ N=0 \} \right] 				\overset{a.s.} \longrightarrow 0. \hspace{15pt} \Box
	\end{equation*}
	
\vspace{15pt}

\textit{Proof of Theorem \ref{Consistency}}\\
\

In what follows, we will be using the well-known inequality 
	
	\begin{equation*}
	L_{n}(\phi_n) - L(\phi_{B}) \leq E \left[ \Big| \eta(\textbf{X}) - 				\eta_{n}(\textbf{X}) \Big| \Bigg| D_n \right] 
	\end{equation*}	
	
to bound the error difference $L_{n}(\phi_n) - L(\phi_{B})$ (see, for example, Devroye et al. (1996)). Here, $\eta(\textbf{X}) = E( Y \vert \textbf{X})$ is given by

\begin{equation*}
\eta(\textbf{X})= \dfrac{\dfrac{p}{\mu(\mathcal{C}_{1})} I\{ \textbf{X} \in \mathcal{C}_{1}\} }{\dfrac{p}{\mu(\mathcal{C}_{1})} I\{ \textbf{X} \in \mathcal{C}_{1}\} +\dfrac{(1-p)}{\mu(\mathcal{C}_{0})} I\{ \textbf{X} \in \mathcal{C}_{0}\}},\\[10pt]
\end{equation*}

whereas the sample-based approximation $\eta_{n}(\textbf{X}, D_{n})$ is  

\begin{equation*}
\eta_{n}(\textbf{X},D_n):=\eta_{n}(\textbf{X})= \dfrac{\dfrac{\widehat{p}}{\mu(\widehat{\mathcal{C}}_{1})} I\{ \textbf{X} \in \widehat{\mathcal{C}}_{1}\} }{\dfrac{\widehat{p}}{\mu(\widehat{\mathcal{C}}_{1})} I\{ \textbf{X} \in \widehat{\mathcal{C}}_{1}\} +\dfrac{(1-\widehat{p})}{\mu(\widehat{\mathcal{C}}_{0})} I\{ \textbf{X} \in \widehat{\mathcal{C}}_{0}\}}.\\[15pt]
\end{equation*}

First, note that \small
	\begin{align}
	\label{Eta - Eta n}
	\Big| \eta(\textbf{X}) - \eta_{n}(\textbf{X}) \Big| &= \Bigg| \dfrac{\dfrac{p}{\mu(\mathcal{C}_{1})} I\{ \textbf{X} \in \mathcal{C}_{1}\} }	{\dfrac{p}{\mu(\mathcal{C}_{1})} I\{ \textbf{X} \in \mathcal{C}_{1}\} +\dfrac{(1-p)}{\mu(\mathcal{C}_{0})} I\{ \textbf{X} \in \mathcal{C}_{0}\}} -\dfrac{\dfrac{\widehat{p}}	{\mu(\widehat{\mathcal{C}}_{1})} I\{ \textbf{X} \in	\widehat{\mathcal{C} 		}_{1}\} }{\dfrac{\widehat{p}}{\mu(\widehat{\mathcal{C}}_{1})} I\{ 				\textbf{X} \in \widehat{\mathcal{C}}_{1}\} +\dfrac{(1-\widehat{p})}				{\mu(\widehat{\mathcal{C}}_{0})} I\{ \textbf{X} \in \widehat{\mathcal{C} 		}_{0}\}} \Bigg| \notag \\[10pt]
	& := \Bigg| \dfrac{A}{B + C} - \dfrac{A'}{B' + C'} \Bigg| \notag \\[10pt]
	& = \Bigg| \dfrac{-A'/(B'+C')}{(B+C)}[(B'+C')-(B+C)] + \dfrac{A'-A}				{(B+C)} 	\Bigg| \notag \\[10pt]
	& \leq c_{0} \Big| A - A' \Big| + c_{1}\Big| C - C' \Big|, 
	\end{align}

\normalsize
where $c_{0}=2 \left[ \min \left\lbrace \frac{p}{\mu(\mathcal{C}_{1})}, \frac{(1-p)}{\mu(\mathcal{C}_{0})} \right\rbrace \right] ^{-1}$ and $c_{1}= \frac{1}{2}c_{0}$. Using a similar manipulation as above, we have 
	\begin{align*}
	\vert A-A' \vert \leq \Bigg| \dfrac{\widehat{p}}{\mu(\widehat{\mathcal{C} 	}_{1})} - \dfrac{p}{\mu(\mathcal{C}_{1})} \Bigg| + \dfrac{1}					{\mu(\widehat{\mathcal{C} }_{1})} \Bigg| I\{ \textbf{X} \in 					\widehat{\mathcal{C}}_{1}\} - I\{ \textbf{X} \in \mathcal{C}_{1}\} 				\Bigg|.
	\end{align*}

But since

	\begin{align*}
	\Big| I\{ \textbf{X} \in \widehat{\mathcal{C}}_{1}\} - I\{ \textbf{X} 			\in \mathcal{C}_{1}\} \Big| &= \Big| I\{ \textbf{X} \in 						\widehat{\mathcal{C}}_{1}\} - I\{ \textbf{X} \in 								[\widehat{\mathcal{C}}_{1} \cup ( \mathcal{C}_{1} \Delta 						\widehat{\mathcal{C}}_{1})]  \} \Big| =  I\{ \textbf{X} \in ( 					\mathcal{C}_{1} \Delta \widehat{\mathcal{C} }_{1})\},
	\end{align*}	

\vspace{7pt}

we find	

	\begin{equation}
	\label{Theorem - A - A'}
	\vert A-A' \vert \leq \Bigg| \dfrac{\widehat{p}}{\mu(\widehat{\mathcal{C} 	}_{1})} - \dfrac{p}{\mu(\mathcal{C}_{1})} \Bigg| + \dfrac{1}					{\mu(\widehat{\mathcal{C} }_{1})} I\{ \textbf{X} \in ( \mathcal{C}_{1} 			\Delta \widehat{\mathcal{C}}_{1})\}.
	\end{equation}
\
	
Similarly, it can be shown that
	\begin{equation}
	\label{Theorem - C - C'}
	\vert C-C' \vert \leq \Bigg| \dfrac{(1-\widehat{p})}{\mu(\widehat{ 				\mathcal{C}}_{0})} - \dfrac{(1-p)}{\mu(\mathcal{C}_{0})} \Bigg| + 				\dfrac{1}{\mu(\widehat{\mathcal{C} }_{0})} I\{ \textbf{X} \in ( 				\mathcal{C}_{0} \Delta \widehat{\mathcal{C}}_{0})\}.
	\end{equation}	
	
Plugging the inequalities in (\ref{Theorem - A - A'}) and (\ref{Theorem - C - C'}) into (\ref{Eta - Eta n}) gives 

	\begin{align*}
	\Big| \eta(\textbf{X}) - \eta_{n}(\textbf{X}) \Big| \leq &c_{0} 				\Bigg| \dfrac{\widehat{p}}{\mu(\widehat{\mathcal{C} 	    					}_{1})} - \dfrac{p}{\mu(\mathcal{C}_{1})} \Bigg| + \dfrac{c_{0}}				{\mu(\widehat{\mathcal{C} }_{1})} I\{ \textbf{X} \in ( \mathcal{C}_{1} 			\Delta \widehat{\mathcal{C}}_{1})\}\\[6pt]
	&+ c_{1}\Bigg| \dfrac{(1-\widehat{p})}{\mu(\widehat{ \mathcal{C}}_{0})} - 	\dfrac{(1-p)}{\mu(\mathcal{C}_{0})} \Bigg| + \dfrac{c_{1}}						{\mu(\widehat{\mathcal{C} }_{0})} I\{ \textbf{X} \in ( 							\mathcal{C}_{0} \Delta \widehat{\mathcal{C}}_{0})\} .
	\end{align*}

Therefore, conditioning on $D_{n}$ and taking expectation, we have
	\begin{align*}
	E \left[ \Big| \eta(\textbf{X}) - \eta_{n}(\textbf{X}) \Big| \Bigg| D_{n} 	\right] & \leq   c_{0} \Bigg| \dfrac{\widehat{p}}{\mu(\widehat{ 				\mathcal{C}}_{1})} - \dfrac{p}{\mu(\mathcal{C}_{1})} \Bigg| + c_{1}  			\Bigg| \dfrac{(1-\widehat{p})}{\mu(\widehat{ \mathcal{C}}_{0})} - 				\dfrac{(1-p)}{\mu(\mathcal{C}_{0})} \Bigg| \\[10pt]
	&+\dfrac{ c_{0}}{\mu(\widehat{\mathcal{C}}_{1})} P \left\lbrace 				\textbf{X} 	\in ( \mathcal{C}_{1} \Delta \widehat{\mathcal{C}}_{1}) 			\Bigg| D_{n} \right\rbrace + \dfrac{c_{1}}{\mu(\widehat{\mathcal{C} 			}_{0})} P \left\lbrace \textbf{X} 	\in ( \mathcal{C}_{0} \Delta 				\widehat{\mathcal{C}}_{0}) \Bigg| D_{n} \right\rbrace\\[10pt]
	&:= A_{1} + A_{2} + A_{3} + A_{4}.
	\end{align*}	
	
Lemma \ref{Lemma - Symmetric Difference}, in conjunction with the strong law of large numbers, implies that both $A_{1}$ and $A_{2}$ converge to $0$ almost surely as $n \to \infty$. Corollary \ref{Corollary - Symmetric Difference} establishes the same result for the terms $A_{3}$ and $A_{4}$, completing the proof. $\hspace{10pt} \Box$

%%%%%%%%%%%%%%%%%%%%%%%%%%%%%%%%%%%%%%%%%%%%%%%%

\section{Classification with Missing Covariates}

%%%%%%%%%%%%%%%%%%%%%%%%%%%%%%%%%%%%%%%%%%%%%%%%

We now consider the more difficult case where some of the covariates can be missing in both the data and observation to be classified. Let $\textbf{Z}=(\textbf{X}',\textbf{V}')' \in \mathbb{R}^{d+s}$ be the vector of covariates to be used to predict the class membership $Y\in \{0,1\}$, where $\textbf{X} \in \mathbb{R}^d, d\geq 1$ is always observable but $\textbf{V} \in \mathbb{R}^s, s \geq 1$ can be missing. Also, let $\delta$ be a $\{0,1\}-$valued random variable defined by $\delta=0$ if $\textbf{V}$ is missing and $\delta=1$ otherwise. \\

We will be assuming that  $\textbf{Z} \vert (\delta=1, Y=1)  \sim $ Unif $(\mathcal{C}_{1})$ and $\textbf{Z} \vert (\delta=1, Y=0)  \sim $ Unif $(\mathcal{C}_{0})$ where $\mathcal{C}_{1}$ and $\mathcal{C}_{0}$ are convex subsets of $\mathbb{R}^{d+s}$. Similarly, let  $\textbf{Z} \vert (\delta=0, Y=1)  \sim $ Unif $(\mathcal{C}^{'}_{1})$ and $\textbf{Z} \vert (\delta=0, Y=0)  \sim $ Unif $(\mathcal{C}^{'}_{0})$ where $\mathcal{C}^{'}_{1}$ and $\mathcal{C}^{'}_{0}$ are the projections of $\mathcal{C}_{1}$ and $\mathcal{C}_{0}$ on $\mathbb{R}^{d}$; that is $\mathcal{C}^{'}_{1} = $ proj$_{\mathbb{R}^d}(\mathcal{C}_{1})$ and $\mathcal{C}^{'}_{0} = $ proj$_{\mathbb{R}^d}(\mathcal{C}_{0})$. Note that both $\mathcal{C}^{'}_{1}$ and $\mathcal{C}^{'}_{0}$ are convex subsets of $\mathbb{R}^{d}$ (see Rockafeller (1996)). \\

Now consider a new observation $(\textbf{Z},Y,\delta)$ that needs to be classified. Mojirsheibani and Montazeri (2007) show that the best classifier (in the sense of having the lowest misclassification error rate) is 

\begin{equation}
\label{Bayes Rule - Missing}
\phi_{B}(\textbf{Z},\delta)= \delta \tilde{\phi}_1(\textbf{Z}) + (1-\delta)\tilde{\phi}_0(\textbf{X}),
\end{equation}

where

\begin{equation}
\label{Bayes missing Piece 1}
\tilde{\phi}_{1}(\textbf{Z})=
\begin{cases}
1 \text{ if } p(\textbf{Z},1)P\{Y=1 \vert \textbf{Z}\} > p(\textbf{Z},0) P\{Y=0 \vert \textbf{Z}\}\\
0 \text{ otherwise, }
\end{cases}
\end{equation}

\begin{equation}
\label{Bayes missing Piece 2}
\tilde{\phi}_{0}(\textbf{X})=
\begin{cases}
1 \text{ if } (1-q(\textbf{X},1))P\{Y=1 \vert \textbf{X}\} > (1 -q(\textbf{X}, 0)) P\{Y=0 \vert \textbf{X}\}\\
0 \text{ otherwise, }
\end{cases}
\end{equation}

\vspace{10pt}

and where 
	\begin{equation*}
	p(\textbf{Z},i)=P\{\delta=1 \vert \textbf{Z}, Y=i\} \hspace{15pt} 		\text{and} \hspace{15pt} q(\textbf{X},i)=P\{\delta=1 \vert 				\textbf{X}, Y=i\}, i=0,1.
	\end{equation*}

The terms $p(\textbf{Z},i)$ and $q(\textbf{X},i)$ are often referred to as the 
missingness probability mechanisms. For more recent results on classification 
with missing covariates, one may also refer to Mojirsheibani (2012) and 
Mojirsheibani and Chenouri (2011).\\
\

To get a more suitable form for the optimal classifier in (\ref{Bayes Rule - Missing}), use the arguments in the proof of Theorem 2 to rewrite $\tilde{\phi}_{1}(\textbf{z})$ in expression (\ref{Bayes missing Piece 1}) as 

\begin{align*}
\tilde{\phi}_{1}(\textbf{z})&= I \left\lbrace \textbf{z} \in \mathcal{C}_{1} \cap \mathcal{C}_{0}, p(\textbf{z},1)P\{Y=1\vert \textbf{Z} = \textbf{z}\} > p(\textbf{z},0) P\{Y=0 \vert \textbf{Z} = \textbf{z}\} \right\rbrace\\
&+ I \left\lbrace \textbf{z} \in \mathcal{C}_{1} - \mathcal{C}_{0}, p(\textbf{z}, 1)>0 \right\rbrace.
\end{align*}

An analogous argument yields

\begin{align*}
\tilde{\phi}_{0}(\textbf{x})&= I \left\lbrace \textbf{x} \in \mathcal{C}^{'}_{1} \cap \mathcal{C}^{'}_{0}, (1-q(\textbf{x},1))P\{Y=1\vert \textbf{X} = \textbf{x}\} > (1-q(\textbf{x},0)) P\{Y=0 \vert \textbf{X} = \textbf{x}\} \right\rbrace\\
&+ I \left\lbrace \textbf{x} \in \mathcal{C}^{'}_{1} - \mathcal{C}^{'}_{0}, q(\textbf{x},1) < 1\right\rbrace.
\end{align*}

Next, let $P\{Y=1\}=p$ and $P\{Y=0\} = 1-p$, where $0 \leq p \leq 1.$ Also recall that $\textbf{Z} \vert (\delta = 1, Y=1) \sim \text{Unif}(\mathcal{C}_{1})$ and $\textbf{Z} \vert (\delta = 1, Y=0) \sim \text{Unif}(\mathcal{C}_{0})$. Using the definition in (\ref{Conditional 1}), we find 
\begin{align}
\label{Bayes Final 1}
\tilde{\phi}_{1}(\textbf{z})&= I \left\lbrace \textbf{z} \in \mathcal{C}_{1} \cap \mathcal{C}_{0}, p(\textbf{z},1) \dfrac{p}{\mu{(\mathcal{C}_{1})}} > p(\textbf{z},0) \dfrac{1-p}{\mu{(\mathcal{C}_{0})}} \right\rbrace \notag \\
&+ I \left\lbrace \textbf{z} \in \mathcal{C}_{1} - \mathcal{C}_{0}, p(\textbf{z},1)>0 \right\rbrace.
\end{align}

A similar approach can be used to find $\tilde{\phi}_{0}(\textbf{X})$ in (\ref{Bayes missing Piece 2}):
\begin{align}
\label{Bayes Final 2}
\tilde{\phi}_{0}(\textbf{x})&= I \left\lbrace \textbf{x} \in \mathcal{C}^{'}_{1} \cap \mathcal{C}^{'}_{0}, (1-q(\textbf{x},1)) \dfrac{p}{\mu{(\mathcal{C}^{'}_{1})}} > (1-q(\textbf{x},0)) \dfrac{1-p}{\mu{(\mathcal{C}^{'}_{0})}} \right\rbrace \notag \\
&+ I \left\lbrace \textbf{x} \in \mathcal{C}^{'}_{1} - \mathcal{C}^{'}_{0}, q(\textbf{x},1) < 1\right\rbrace.
\end{align}

Combining the classifiers in (\ref{Bayes Final 1}) and (\ref{Bayes Final 2}) gives an expression for the optimal classifier, as expressed in the following theorem.

\begin{theorem}
Let $\textbf{Z} \vert (\delta=1, Y=1)  \sim $ Unif $(\mathcal{C}_{1})$ and $\textbf{Z} \vert (\delta=1, Y=0)  \sim $ Unif $(\mathcal{C}_{0})$ where $\mathcal{C}_{1}$ and $\mathcal{C}_{0}$ are convex subsets of $\mathbb{R}^{d+s}$. Similarly, let  $\textbf{Z} \vert (\delta=0, Y=1)  \sim  $ Unif $ (\mathcal{C}^{'}_{1})$ and $\textbf{Z} \vert (\delta=0, Y=0)  \sim  $ Unif $ (\mathcal{C}^{'}_{0})$ where $\mathcal{C}^{'}_{1} = $ proj$_{\mathbb{R}^d}(\mathcal{C}_{1})$ and $\mathcal{C}^{'}_{0} = $ proj$_{\mathbb{R}^d}(\mathcal{C}_{0})$. Also, assume that 
\begin{equation*}
\inf_{\textbf{z} \in \mathbb{R}^{d+s}} P \{ \delta = 1 \vert \textbf{Z} = \textbf{z}, Y=1\} > 0 \hspace{12pt} \text{and} \hspace{12pt}  \sup_{\textbf{x} \in \mathbb{R}^{d}} P \{ \delta = 1 \vert \textbf{X} = \textbf{x}, Y=1\} < 1.
\end{equation*}

Then the optimal classifier is given by

\begin{equation}
\label{Bayes Missing Convex}
\phi_{B}(\textbf{Z},\delta)=\delta \tilde{\phi}_{1}(\textbf{z}) +(1-\delta)\tilde{\phi}_{0}(\textbf{x})\\[10pt]
\end{equation}

where
\begin{align*}
\tilde{\phi}_{1}(\textbf{z})= I \left\lbrace \textbf{z} \in \mathcal{C}_{1} \cap \mathcal{C}_{0}, p(\textbf{z},1) \dfrac{p}{\mu{(\mathcal{C}_{1})}} > p(\textbf{z},0) \dfrac{1-p}{\mu{(\mathcal{C}_{0})}} \right\rbrace + I \left\lbrace \textbf{z} \in \mathcal{C}_{1} - \mathcal{C}_{0} \right\rbrace
\end{align*}

and
\begin{align*}
\tilde{\phi}_{0}(\textbf{x})= I \left\lbrace \textbf{x} \in \mathcal{C}^{'}_{1} \cap \mathcal{C}^{'}_{0}, (1-q(\textbf{x},1)) \dfrac{p}{\mu{(\mathcal{C}^{'}_{1})}} > (1-q(\textbf{x},0)) \dfrac{1-p}{\mu{(\mathcal{C}^{'}_{0})}} \right\rbrace  + I \left\lbrace \textbf{x} \in \mathcal{C}^{'}_{1} - \mathcal{C}^{'}_{0}\right\rbrace.
\end{align*}
\end{theorem}

\vspace{10pt}

Since it is unrealistic to assume that the four convex sets and the functional forms of the missingness probability mechanisms $p(\textbf{z},y)$ and $q(\textbf{x},y)$ are known, we must proceed by finding a suitable plug-in version of the rule in (\ref{Bayes Missing Convex}). Let $\widehat{\mathcal{C}}_{1}$, $\widehat{\mathcal{C}}_{0}$, $\widehat{\mathcal{C}}_{1} \hspace{.5pt}^{'}$ and $\widehat{\mathcal{C}}_{0} \hspace{.5pt}^{'}$ be the sample based counterparts of $\mathcal{C}_{1}$, $\mathcal{C}_{0}$, $\mathcal{C}^{'}_{1}$ and $\mathcal{C}^{'}_{0}$ respectively; that is, let
	\begin{align*}
	\widehat{\mathcal{C}}_{1} &= \text{ convex hull of } \{(\textbf{X}_{i}, 		\textbf{V}_{i}) \vert \delta_{i}=1, Y_{i} =1 \}_{i=1}^{n},\\[5pt]
	\widehat{\mathcal{C}}_{0} &= \text{ convex hull of } \{(\textbf{X}_{i}, 		\textbf{V}_{i}) \vert \delta_{i}=1, Y_{i} =0 \}_{i=1}^{n},\\[5pt]
	\widehat{\mathcal{C}}_{1} \hspace{.5pt}^{'} &= \text{ convex hull of } 			\{\textbf{X}_{i} \vert \delta_{i}=0, Y_{i} =1 \}_{i=1}^{n},\\[5pt]
	\widehat{\mathcal{C}}_{0} \hspace{.5pt}^{'} &= \text{ convex hull of } 			\{\textbf{X}_{i} \vert \delta_{i}=0, Y_{i} =0 \}_{i=1}^{n}.
	\end{align*}

We propose a data-based version of the classifier in (\ref{Bayes Missing Convex}) given by

	\begin{equation}
	\label{Naive Bayes Rule Data Version}
	\phi_{n}(\textbf{Z},\delta)= \delta \phi_{n,1}(\textbf{Z}) + 					(1-\delta)\phi_{n,0}(\textbf{X})
	\end{equation}

where

	\begin{equation}
	\label{Bayes Convex Sets No Missing 2}
	\phi_{n,1}(\textbf{z})=
		\begin{cases}
		1 \text{ if } \textbf{z} \in \widehat{\mathcal{C}}_{1} - \widehat 				{\mathcal{C}}_{0}\\
		0 \text{ if } \textbf{z} \in \widehat{\mathcal{C}}_{0} - \widehat{ 				\mathcal{C}}_{1}\\
		1 \text{ if } \textbf{z} \in \widehat{\mathcal{C}}_{1}	\cap 					\widehat{\mathcal{C}}_{0}, \hspace{5pt}	\widehat{p}(\textbf{z},1)				\dfrac{ \widehat{p}}{\mu(\widehat{\mathcal{C}}_{1})} > \widehat{p} 				 (\textbf{z},0) \dfrac{1-\widehat{p}}{\mu(\widehat{ \mathcal{C} 				}_{0})}\\[15pt]
		0 \text{ if } \textbf{z} \in \widehat{\mathcal{C}}_{1}	\cap 					\widehat{\mathcal{C}}_{0}, \hspace{5pt}	\widehat{p}(\textbf{z},1) 				\dfrac{ \widehat{p}}{\mu(\widehat{\mathcal{C}}_{1})} \leq \widehat{p} 		(\textbf{z},0) \dfrac{1-\widehat{p}}{\mu(\widehat{ \mathcal{C} 				    }_{0})}\\[15pt]
		I\left\lbrace \widehat{p} \cdot \inf\limits_{\textbf{y}\in 						\widehat{\mathcal{C}}_{1}}\| \textbf{z}-\textbf{y}\| < (1 - 					\widehat{p}) \cdot \inf\limits_{\textbf{y}\in	\widehat{\mathcal{C} 			}_{0}}\|\textbf{z}-\textbf{y}\| \right\rbrace \hspace{15pt} \text{ if 		} \hspace{5pt} \textbf{z} \not\in \widehat{\mathcal{C}}_{1} \cup 				\widehat{\mathcal{C}}_{0},
		\end{cases}
	\end{equation}

and

	\begin{equation}
	\label{Bayes Convex Sets Missing Data 2} 
	\phi_{n,0}(\textbf{x})=
		\begin{cases}
		1 \text{ if } \textbf{x} \in \widehat{\mathcal{C}}_{1} 							\hspace{.5pt}^{'} - \widehat{\mathcal{C}}_{0} \hspace{.5pt}^{'}\\
		0 \text{ if } \textbf{x} \in \widehat{\mathcal{C}}_{0} 							\hspace{.5pt}^{'} - \widehat{\mathcal{C}}_{1} \hspace{.5pt}^{'}\\
		1 \text{ if } \textbf{x} \in \widehat{\mathcal{C}}_{1} 							\hspace{.5pt}^{'} \cap \widehat{\mathcal{C}}_{0}  \hspace{.5pt}^{'}, 			\left( 1 - \widehat{q} (\textbf{x},1)\right) \dfrac{\widehat{p}}				{\mu(\widehat{\mathcal{C} }_{1} \hspace{.5pt}^{'})} > \left( 					1-\widehat{q}(\textbf{x},0) \right)\dfrac{1-\widehat{p}}						{\mu(\widehat{\mathcal{C}}_{0} \hspace{.5pt}^{'} )}\\[15pt]
		0 \text{ if } \textbf{x} \in \widehat{\mathcal{C}}_{1} 							\hspace{.5pt}^{'} \cap \widehat{\mathcal{C}}_{0}  \hspace{.5pt}^{'}, 			\left( 1 - \widehat{q} (\textbf{x},1)\right) \dfrac{\widehat{p}}				{\mu(\widehat{\mathcal{C} }_{1} \hspace{.5pt}^{'})} \leq \left( 				1-\widehat{q}(\textbf{x},0) \right)\dfrac{1-\widehat{p}}						{\mu(\widehat{\mathcal{C}}_{0} \hspace{.5pt}^{'} )}\\[15pt]
		I\left\lbrace \widehat{p} \cdot \inf\limits_{\textbf{y}\in 						\widehat{\mathcal{C}}_{1} \hspace{.5pt}^{'}} \|									\textbf{x}-\textbf{y}\| < (1 - \widehat{p}) \cdot 								\inf\limits_{\textbf{y}\in \widehat{\mathcal{C}}_{0} 							\hspace{.5pt}^{'}} \| \textbf{x}-\textbf{y}\| \right\rbrace 					\hspace{15pt} \text{ if } \hspace{5pt} \textbf{x} \not\in 						\widehat{\mathcal{C}}_{1} \hspace{.5pt}^{'} \cup 								\widehat{\mathcal{C}}_{0} \hspace{.5pt}^{'}.
		\end{cases}
	\end{equation}

Here, $\widehat{p}={n}^{-1}\sum_{i=1}^{n} Y_{i}$. As for the terms $\widehat{p}(\textbf{z},y)$ and $\widehat{q}( \textbf{x},y)$, we consider a  method of estimation based on kernel regression.\\
\

Recall the forms of the missingness probability mechanisms $p_{y}(\textbf{z}) := p(\textbf{z},y)$ and $q_{y}(\textbf{x}) := q(\textbf{x},y)$:
	\begin{align*}
	p_{y}(\textbf{z}) &= P\{\delta = 1 \vert \textbf{Z} = \textbf{z}, Y=y\} 		=E(\delta \vert \textbf{Z}=\textbf{z},Y=y),\\[5pt]
	q_{y}(\textbf{x}) &= P\{\delta = 1 \vert \textbf{X} = \textbf{x}, Y=y\} 		=E(\delta \vert \textbf{X}=\textbf{x},Y=y).
	\end{align*}

Under the commonly used assumption of data missing at random (MAR), one assumes that the probability that $\textbf{V}$ is missing does not depend on $\textbf{V}$ itself. That is,
	\begin{equation*}
	p_{y}(\textbf{z}) = P\{\delta = 1 \vert \textbf{Z}=\textbf{z},Y=y\} = 			P\{\delta = 1 \vert \textbf{X}=\textbf{x},Y=y\}= q_{y}(\textbf{x}).
	\end{equation*}

Under the MAR assumption, the kernel regression estimates of $q_{y}(\textbf{x})$, $y=0,1,$ are

	\begin{equation}
	\label{Kernel Regression under MAR}
	\widehat{q}_{y}(\textbf{x}) = \sum_{j:Y_j=y} \delta_{j} K\left( 				\dfrac{\textbf{X}_j - \textbf{x}}{h_{n}} \right) \Big/ 							\sum_{j:Y_j=y}  K\left( \dfrac{\textbf{X}_j - \textbf{x}}{h_{n}} \right),
	\end{equation}

where $0/0$ is defined to be $0$. Here, the kernel $K$ is a map of the form $K:\mathbb{R}^{d} \to \mathbb{R}^{+}$ with  smoothing parameter $h_{n}$ ($h_{n} \to 0$ as $n \to \infty$).\\
\

To evaluate the performance of the rule $\phi_{n}$ in (\ref{Naive Bayes Rule Data Version}), let
	\begin{equation*}
	L(\phi_{B}) = P \{ \phi_{B} (\textbf{Z},\delta) \neq Y\} \hspace{15pt} 			\text{ and } \hspace{15pt} L_{n}(\phi_{n}) = P \{ \phi_{n}(\textbf{Z},			\delta) \neq Y \vert D_{n} \}.
	\end{equation*}

In what follows, we shall assume that the chosen kernel in (\ref{Kernel Regression under MAR}) is \textit{regular:} a nonnegative kernel $K$ is said to be regular if there are positive constants $b>0$ and $r>0$ for which $K(\textbf{x})\geq bI\{ \textbf{x}\in B(0,r)\}$ and $\int \sup_{\textbf{y}\in \textbf{x}+ B(0,r)}K(\textbf{y})d\textbf{x}<\infty,$ where $B(0,r)$ is the ball of radius $r$ centered at the origin (for more on this see, for example, Gy\"{o}rfi et al. (2002)).

\begin{theorem}
\label{Missing Consistency}
Let $\phi_{n}$ be the classifier in (\ref{Naive Bayes Rule Data Version}) where the kernel $K$ is regular and the bandwidth $h_{n}$ satisfies $h_{n} \to 0$ and $nh_{n}^{d} \to \infty$. Then $L_{n}(\phi_{n}) \overset{a.s.} \longrightarrow L(\phi_{B})$ as $n \to \infty$, where $\phi_{B}$ is the optimal rule in (\ref{Bayes Missing Convex}).
\end{theorem}

To prove this result, we need the following lemma.

\begin{lemma}[Mojirsheibani and Montazeri (2007)]
\label{Lemma - Missing}
Let $\phi_{B}$ be the optimal rule given by (\ref{Bayes Missing Convex}). For $j=0,1$, let $f_{j}:\mathbb{R}^{d+s} \to [0,1]$ and $g_{j}: \mathbb{R}^{d} \to [0,1]$ be any given functions. Furthermore, let $\phi_{1}(\textbf{z})=I\{f_{1}(\textbf{z})>f_{0}(\textbf{z})\}$ and $\phi_{0}(\textbf{x})=I\{g_{1}(\textbf{x})>g_{0}(\textbf{x})\}$. Also, put
	\begin{equation*}
	\phi(\textbf{Z},Y)=\delta \phi_{1}(\textbf{Z}) + (1-\delta) \phi_{0}			(\textbf{X}).
	\end{equation*}
Then
	\begin{align*}
	L(\phi) - L(\phi_{B}) &\leq 2 \sum_{i=0}^{1} E \left[ \Big| 					p(\textbf{Z},i)P(Y=i \vert \textbf{Z}) - f_{i}(\textbf{Z})  \Big| 				\right]\\
	 &+ 2 \sum_{i=0}^{1} E \left[ \Big| (1-q(\textbf{X},i))P(Y=i \vert 				\textbf{X}) - g_{i}(\textbf{X})  \Big| \right].
	\end{align*}
\end{lemma}

\textit{Proof of Theorem \ref{Missing Consistency}}\\
\

We begin by introducing some notation. Let

\begin{equation*}
f_{1}(\textbf{Z})=\widehat{p}(\textbf{Z},1)  					\dfrac{\dfrac{\widehat{p}}{\mu(\widehat{\mathcal{C}}_{1})} I\{ \textbf{Z} \in \widehat{\mathcal{C}}_{1}\} }{\dfrac{\widehat{p}}{\mu(\widehat{\mathcal{C}}_{1})} I\{ \textbf{Z} \in \widehat{\mathcal{C}}_{1}\} +\dfrac{(1-\widehat{p})}{\mu(\widehat{\mathcal{C}}_{0})} I\{ \textbf{Z} \in \widehat{\mathcal{C}}_{0}\}}, \\[5pt]
\end{equation*}
\begin{equation*}
f_{0}(\textbf{Z})=\widehat{p}(\textbf{Z},0)  					\dfrac{\dfrac{(1-\widehat{p})}{\mu(\widehat{\mathcal{C}}_{0})} I\{ \textbf{Z} \in \widehat{\mathcal{C}}_{0}\} }{\dfrac{\widehat{p}}{\mu(\widehat{\mathcal{C}}_{1})} I\{ \textbf{Z} \in \widehat{\mathcal{C}}_{1}\} +\dfrac{(1-\widehat{p})}{\mu(\widehat{\mathcal{C}}_{0})} I\{ \textbf{Z} \in \widehat{\mathcal{C}}_{0}\}}, \\[5pt]
\end{equation*}
\begin{equation*}
g_{1}(\textbf{X})=(1-\widehat{q}(\textbf{X},1))  					\dfrac{\dfrac{\widehat{p}}{\mu(\widehat{\mathcal{C}}_{1} \hspace{.5pt}^{'})} I\{ \textbf{X} \in \widehat{\mathcal{C}}_{1} \hspace{.5pt}^{'}\} }{\dfrac{\widehat{p}}{\mu(\widehat{\mathcal{C}}_{1} \hspace{.5pt}^{'})} I\{ \textbf{X} \in \widehat{\mathcal{C}}_{1} \hspace{.5pt}^{'}\} +\dfrac{(1-\widehat{p})}{\mu(\widehat{\mathcal{C}}_{0} \hspace{.5pt}^{'})} I\{ \textbf{X} \in \widehat{\mathcal{C}}_{0} \hspace{.5pt}^{'}\}}, \\[5pt]
\end{equation*}
\begin{equation*}
g_{0}(\textbf{X})=(1-\widehat{q}(\textbf{X},0))  					\dfrac{\dfrac{(1-\widehat{p})}{\mu(\widehat{\mathcal{C}}_{0} \hspace{.5pt}^{'})} I\{ \textbf{X} \in \widehat{\mathcal{C}}_{0} \hspace{.5pt}^{'}\} }{\dfrac{\widehat{p}}{\mu(\widehat{\mathcal{C}}_{1} \hspace{.5pt}^{'})} I\{ \textbf{X} \in \widehat{\mathcal{C}}_{1} \hspace{.5pt}^{'}\} +\dfrac{(1-\widehat{p})}{\mu(\widehat{\mathcal{C}}_{0} \hspace{.5pt}^{'})} I\{ \textbf{X} \in \widehat{\mathcal{C}}_{0} \hspace{.5pt}^{'}\}}. 
\end{equation*}

\vspace{15pt}
A simple extension of Lemma \ref{Lemma - Missing} yields

	\begin{align}
	\label{Missing Bound}
	L_{n}(\phi_{n}) - L(\phi_{B}) &\leq 2 \sum_{i=0}^{1} E \left[ \Big| 			p(\textbf{Z},i)P(Y=i \vert \textbf{Z}) - f_{i}(\textbf{Z})  \Big| 				\Bigg| D_{n}  \right] \notag \\
	 &+ 2 \sum_{i=0}^{1} E \left[ \Big| (1-q(\textbf{X},i))P(Y=i \vert 				\textbf{X}) - g_{i}(\textbf{X})  \Big| \Bigg| D_{n} \right].
	\end{align}

To prove the theorem, we will show that 

	\begin{equation*}
	E \left[ \Big| 	p(\textbf{Z},1)P(Y=1 \vert \textbf{Z}) - f_{1}					(\textbf{Z}) \Big| \Bigg| D_{n} \right] \overset{a.s.} \longrightarrow 0		\end{equation*}	
	
as $n \to \infty$; the other terms in (\ref{Missing Bound}) can be dealt with analogously. Following the proof of Theorem \ref{Consistency}, we find
	\begin{equation*}
	\label{Big inequality}
	\Big| p(\textbf{Z},1)P(Y=1 \vert \textbf{Z}) - f_{1}(\textbf{Z})\Big| \leq c_{0} \left\lbrace \Big| B - B' \Big| + \Big| C - C' \Big| + \Big| D 	- D' \Big| \right\rbrace, \\[10pt] 
	\end{equation*}

where 
	\begin{equation*}
	c_{0}= \left[ \min \left\lbrace \frac{p}{\mu(\mathcal{C}_{1})}, \frac{(1-p)}{\mu(\mathcal{C}_{0})} \right\rbrace \right] ^{-1}, \\[10pt]
	\end{equation*}
and where $B, B', C,$ and $C'$ are defined as before, with $\textbf{Z}$ playing the role of $\textbf{X}$. Here, as before,

	\begin{equation*}
	\Big| B - B' \Big| \leq \Bigg| \frac{\widehat{p}}{\mu(\widehat{ 				\mathcal{C} }_{1})} - \frac{p}{\mu(\mathcal{C}_{1})} \Bigg| + \frac{1}			{\mu(\widehat{ \mathcal{C}}_{1})} I \{ \textbf{Z} \in (\mathcal{C}_{1} 			\Delta \widehat{\mathcal{C}}_{1})\}
	\end{equation*}

and
	\begin{equation*}
	\Big| C-C' \Big| \leq \Bigg| \dfrac{(1-\widehat{p})}{\mu(\widehat{ \mathcal{C}}_{0})} - \dfrac{(1-p)}{\mu(\mathcal{C}_{0})} \Bigg| + 				\dfrac{1}{\mu(\widehat{\mathcal{C} }_{0})} I\{ \textbf{Z} \in (\mathcal{C}_{0} \Delta \widehat{\mathcal{C}}_{0})\}. \\[10pt] 
	\end{equation*}

As for the term $ \vert D - D' \vert$, it is easy (but tedious) to show that

	\begin{align*}
	\Big| D - D' \Big| &= \Bigg| \frac{p(\textbf{Z},1)}{ [\mu(\mathcal{ C}_{1}) / p \cdot I \{ \textbf{Z} \in \mathcal{C}_{1}  \} ] } - 				\frac{\widehat{p}(\textbf{Z},1)}{ [\mu( \widehat{\mathcal{C}}_{1}) / \widehat{p} \cdot I \{ \textbf{Z} \in \widehat{\mathcal{C}}_{1}  \} ] } 		\Bigg| \\[7pt]
	& \leq \Bigg| \frac{\widehat{p}}{\mu(\widehat{ \mathcal{C} }_{1})} - 			\frac{p}{\mu(\mathcal{C}_{1})} \Bigg| + \frac{1} {\mu(\widehat{ 				\mathcal{C}}_{1})} I \{ \textbf{Z} \in (\mathcal{C}_{1} \Delta 					\widehat{\mathcal{C}}_{1})\} + \frac{1}{\mu(\mathcal{C}_{1})} \Big| 			\widehat{p}(\textbf{Z},1) - p(\textbf{Z},1) \Big|. 
	\end{align*}

It follows by (\ref{Big inequality}) that
	\begin{align*}
	\Big| p(\textbf{Z},1)P(Y=1 \vert \textbf{Z}) - f_{1}(\textbf{Z})\Big| 			\leq &2c_{0} \Bigg| \dfrac{\widehat{p}}{\mu(\widehat{\mathcal{C} 	    		}_{1})} - \dfrac{p}{\mu(\mathcal{C}_{1})} \Bigg| + \dfrac{2c_{0}}				{\mu(\widehat{\mathcal{C} }_{1})} I\{ \textbf{Z} \in ( \mathcal{C}_{1} 			\Delta \widehat{\mathcal{C}}_{1})\}\\[6pt]
	&+ c_{0}\Bigg| \dfrac{(1-\widehat{p})}{\mu(\widehat{ \mathcal{C}}_{0})} - 		\dfrac{(1-p)}{\mu(\mathcal{C}_{0})} \Bigg| + \dfrac{c_{0}}						{\mu(\widehat{\mathcal{C} }_{0})} I\{ \textbf{Z} \in ( \mathcal{C}_{0} \Delta \widehat{\mathcal{C}}_{0})\}\\[6pt]
	& +  \frac{c_{0}}{\mu(\mathcal{C}_{1})} \Big| \widehat{p}(\textbf{Z},1) 		- 	p(\textbf{Z},1) \Big|.
	\end{align*}	

Therefore, conditioning on $D_{n}$ and taking expectation, we have

	\begin{align*}
	E \left[ \Big| p(\textbf{Z},1)P(Y=1 \vert \textbf{Z}) - f_{1}					(\textbf{Z})\Big| \Bigg| D_{n} \right] & \leq   2c_{0} \Bigg| 					\dfrac{\widehat{p}}{\mu(\widehat{ \mathcal{C}}_{1})} - \dfrac{p}				{\mu(\mathcal{C}_{1})} \Bigg| + c_{0}  \Bigg| \dfrac{(1-\widehat{p})}			{\mu(\widehat{ \mathcal{C}}_{0})} - \dfrac{(1-p)}{\mu(\mathcal{C}_{0})} 		\Bigg| \\[10pt]
	&+\dfrac{ 2c_{0}}{\mu(\widehat{\mathcal{C}}_{1})} P \left\lbrace 				\textbf{Z} 	\in ( \mathcal{C}_{1} \Delta \widehat{\mathcal{C}}_{1}) 			\Bigg| D_{n} \right\rbrace\\[10pt]
	&+ \dfrac{c_{0}}{\mu(\widehat{\mathcal{C} }_{0})} P \left\lbrace 				\textbf{Z} 	\in ( \mathcal{C}_{0} \Delta \widehat{\mathcal{C}}_{0}) \Bigg| 	D_{n} \right\rbrace\\[10pt]
	& +  \frac{c_{0}}{\mu(\mathcal{C}_{1})} E \left[ \Big| \widehat{p}				(\textbf{Z},1) - p(\textbf{Z},1) \Big| \Bigg| D_{n} \right].
	\end{align*}	
The strong consistency of the kernel regression estimate $\widehat{p}(\textbf{Z},1)$ implies that 

\begin{equation*}
E \left[ \Big| \widehat{p}(\textbf{Z},1) - p(\textbf{Z},1) \Big| \Bigg| D_{n} \right] \overset{a.s.} \longrightarrow 0
\end{equation*}

as $n \to \infty$ (for more on the properties of kernel regression estimates, see, for example, Gy\"{o}rfi et al. (2002)); the remaining terms were dealt with in the proof of Theorem \ref{Consistency}. \hspace{15pt} $\Box$

\vspace{25pt}

\textbf{Acknowledgment}\\
\

This research was supported by a grant to L. Demirdjian from the Interdisciplinary Research Institute of the Sciences at California State University Northridge.

\end{document}